\documentclass{amsart}

\title[Simplicial homology in Mathematica]{SimplicialHomology: implementation of abstract simplicial complex in Mathematica}
\author[Taggar]{Naman Taggar}
\address{Department of Mathematics, University of Delhi, Delhi 110007, India}
\email{namtgr@gmail.com, namtgr@maths.du.ac.in}
\date{14 August, 2026}

\usepackage{tikz}
\usetikzlibrary{positioning,arrows.meta}
\usepackage{url,hyperref}

\usepackage{listings}
\usepackage{xcolor}
\usepackage{lstautogobble}
\usepackage{array,booktabs}

\begin{document}

	\maketitle
	
	\begin{abstract}
		Computations in simplicial homology are crucial for algebraic topology and computational geometry in general. Although several mature software systems provide facilities for working with simplicial complexes, there is currently limited support for such computations within the Wolfram language ecosystem. This paper presents SimplicialHomology, an open-source Wolfram language paclet for constructing abstract simplicial complexes and computing their homology groups over the integers, rational numbers, and finite fields. Furthermore, the package supports many other functionalities such as calculating Betti numbers, automorphism groups, and other constructions such as joins, stars, links, etc. Several examples are presented to illustrate the package and its performance.
		\paragraph{Keywords} Algebraic topology, simplicial complexes, simplicial homology, homology theory, computational topology, Wolfram language.
		\paragraph{AMS Classification} 55N35, 57Q05, 68W30.
	\end{abstract}
	
	\tableofcontents

	\section{Introduction}

	Simplicial complexes provide a combinatorial model for studying topological spaces and play a central role in algebraic topology \cite{Hatcher}. Their discrete nature makes them well suited for computation, and algorithms based on simplicial complexes appear in diverse areas including computational topology, geometric modeling, and topological data analysis. Among the most fundamental invariants associated with a simplicial complex are its homology groups, which capture information about connected components, holes, and higher-dimensional voids in the associated space. A number of software systems support computations involving simplicial complexes. General-purpose mathematical softwares such as SageMath \cite{sagemath} include extensive functionality for simplicial complexes and homology computations, and specialised packages such as GAP \cite{GAP} with HAP \cite{HAP}, GUDHI \cite{gudhi}, and polymake \cite{polymake} provide tools for computational topology and related applications in particular. Despite the extensive capabilities of Mathematica \cite{Mathematica} in symbolic computation, graph theory, linear algebra, and computational geometry, there is comparatively very little in-built or third party support for simplicial homology or discrete topology in general. Existing functionality largely focuses on geometric meshes rather than abstract simplicial complexes and their associated algebraic invariants.
	
	The aim of the SimplicialHomology paclet is to provide a lightweight and idiomatic implementation of simplicial homology for the Wolfram Language. The package supports construction of simplicial complexes from maximal simplices, mesh objects, and built-in triangulations, together with common operations such as joins, cones, suspensions, stars, links, etc. Homology groups may be computed over the integers, rational numbers, or finite fields, and related invariants including Betti numbers and Euler characteristics are also available. The implementation emphasises both usability and efficiency. Boundary operators are stored as sparse matrices, repeated computations are avoided through caching mechanisms, and integer homology is computed using Smith normal form to recover both free and torsion components. These design choices allow the package to handle many examples encountered in teaching and research while remaining entirely within the Wolfram Language.
	
	The remainder of the paper is organised as follows. Section \ref{mb} introduces the necessary mathematical background. Section \ref{des} describes the design and implementation of the package. Section \ref{ex} presents examples and performance benchmarks, and Section \ref{con} concludes the paper. Appendix \ref{a} provides information on the availability of the package.
	
	\section{Mathematical Background}
	\label{mb}
	
	An abstract simplicial complex is a finite collection $K$ of finite subsets of a vertex set $V$ such that whenever $\sigma \in K$ and $\tau \subseteq \sigma$, then $\tau \in K$. The elements of $K$ are called simplices, and a simplex containing $p+1$ vertices is said to have dimension $p$. The maximal simplices of $K$, referred to as facets, uniquely determine the complex and therefore provide a convenient representation for computation.
	
	For each integer $p \geq 0$, let $C_p(K)$ denote the free abelian group generated by the oriented $p$-simplices of $K$. The boundary homomorphism
	$$\partial_p : C_p(K) \longrightarrow C_{p-1}(K)$$
	is defined by
	$$[v_0,\ldots,v_p]\mapsto\sum_{i=0}^{p} (-1)^i [v_0,\ldots,\widehat{v_i},\ldots,v_p],$$
	where $\widehat{v_i}$ indicates that the vertex $v_i$ is omitted. Extending this map linearly yields the chain complex
	$$\cdots\longrightarrow C_{p+1}(K)\xrightarrow{\partial_{p+1}}C_p(K)\xrightarrow{\partial_p}C_{p-1}(K)\longrightarrow\cdots,$$
	which satisfies the fundamental identity
	$$\partial_p\partial_{p+1}=0.$$
	The $p$-th homology group of $K$ is defined by
	$$H_p(K) = \ker(\partial_p)\big/\operatorname{im}(\partial_{p+1}),$$
	and measures the $p$-dimensional topological features of the complex. When coefficients are taken over the integers, the homology groups decompose into free and torsion components. The rank of the free part is the $p$-th Betti number, while the torsion coefficients encode finite-order topological information.
	
	By the fundamental theorem of finitely generated abelian groups, the integer homology groups decompose as
	$$H_p(K; \mathbb{Z}) \cong \mathbb{Z}^{\beta_p} \oplus \mathbb{Z}_{q_1} \oplus \cdots \oplus \mathbb{Z}_{q_k},$$
	where $\beta_p$ is the $p$-th Betti number (the rank of the free part) and the integers $q_i \geq 2$ are the torsion coefficients which encode finite-order topological information. To transition from abstract algebra to computation, we impose a strict global ordering on the vertex set $V$ inherited from Mathematica's usual handling, allowing us to uniquely orient each simplex $[v_0, v_1, \ldots, v_p]$ such that $v_0 < v_1 < \ldots < v_p$. By fixing an ordering of the $p$-simplices, they form a canonical basis for the free abelian group $C_p(K) \cong \mathbb{Z}^{n_p}$, where $n_p$ is the number of $p$-simplices in $K$. With these bases chosen, the boundary homomorphism $\partial_p$ can be explicitly represented as an $n_{p-1} \times n_p$ integer matrix $M_p$. The entry $(M_p)_{i,j}$ is $+1$, $-1$, or $0$, depending on whether the $i$-th $(p-1)$-simplex is a face of the $j$-th $p$-simplex and its corresponding orientation sign. The chain complex identity $\partial_p\partial_{p+1}=0$ translates to the matrix product $M_p M_{p+1} = 0$. 
	
	Computing the homology group $H_p(K)$ over the integers then reduces to standard matrix operations. By applying integer row and column operations to compute the Smith normal form (SNF) of $M_{p+1}$, we obtain a diagonal matrix whose non-zero entries $q_i \geq 2$ yield the torsion coefficients of $H_p(K)$. The $p$-th Betti number $\beta_p$ is efficiently computed via the rank-nullity theorem:
	$$ \beta_p = \dim(\ker M_p) - \operatorname{rank}(M_{p+1}) = (n_p - \operatorname{rank}(M_p)) - \operatorname{rank}(M_{p+1}). $$
	
	The package presented in this paper computes homology over the integers, the rational numbers, and finite fields of prime order. Integer homology is obtained using the Smith normal form of the boundary matrices, allowing both the free and torsion components of the homology groups to be recovered.
	
	\section{Package Design and Implementation}
	\label{des}
	
	SimplicialHomology is implemented entirely in the Wolfram Language as a paclet and follows the standard paclet architecture. The source code is organised into independent modules responsible for constructing simplicial complexes, computing homological invariants, providing utility functions, etc. The package provides built-in triangulations of various spaces such as projective planes \cite{complex}, the Poincar\'e homology $3$-sphere due to Bj\"orner and Lutz \cite{poincare}, the Rudin ball \cite{rudin}, the Ziegler ball \cite{ziegler}, the dunce hat \cite{dunce}, the Cs\'asz\'ar's torus \cite{lutz}, and more. The principal object of the package is {\ttfamily SimplicialComplex}, which represents an abstract simplicial complex together with metadata required for efficient computation. A complex may be constructed from a collection of facets, a {\ttfamily MeshRegion} or {\ttfamily BoundaryMeshRegion}, or one of the built-in triangulations provided by the package. During construction, the simplicial closure is generated automatically. Homology computations are performed using the chain complex associated with the simplicial complex. Boundary operators are represented as sparse matrices, substantially reducing memory usage for complexes containing large numbers of simplices while also improving the efficiency of matrix operations. Over the integers, homology groups are computed using the Smith normal form of the boundary matrices, from which both the free and torsion components are obtained. Computations over the rational numbers and finite fields are performed using the corresponding coefficient arithmetic.

	To improve performance, the package employs an internal caching mechanism. Every simplicial complex is assigned a unique identifier, and computationally expensive quantities such as boundary matrices, matrix ranks, Smith normal forms, homology groups, and Betti numbers are computed only once. Subsequent evaluations retrieve these values directly from the cache, significantly reducing computation time when several invariants of the same complex are requested. In addition to homology computations, the package provides operations including simplicial joins, cones, suspensions, connected components, Euler characteristics, and various combinatorial properties of simplicial complexes. The implementation emphasises compatibility with the symbolic programming nature of the Wolfram Language.
	
	\section{Examples and usage}
	\label{ex}
	
	An abstract simplicial complex is defined directly using the facets constituting it. Ordinarily, maximal facets should be supplied, and the associated closure is automatically generated. The examples in this section illustrate how to do this, compute their invariants, and apply topological operations. The package must be installed and imported by running
	\begin{lstlisting}
		In[1]:= PacletInstall["Taggar/SimplicialHomology"];
		In[2]:= Needs["Taggar`SimplicialHomology`"];
	\end{lstlisting}
	
	\subsection{Basic operations}
	
	For example, I construct a simple triangulation of $S^1$ as a hollow triangle:
	\begin{lstlisting}
		In[3]:= s1 = SimplicialComplex[{{1, 2}, {2, 3}, {1, 3}}];
	\end{lstlisting}
	Once an abstract simplicial complex is defined, the various functionalities offered by SimplicialHomology paclet can be used on this object. Below, we determine the Euler characteristic and $f$-vector of $S^1$:
	\begin{lstlisting}
		In[4]:= s1["EulerCharacteristic"]
		Out[4]= 0
		In[5]:= s1["FVector"]
		Out[5]= {1, 3, 3}
	\end{lstlisting}
	Alternatively, users can call well-known spaces by their name specification due to the package's built-in library of complexes. For example, it provides the minimal triangulation of the 2-torus \cite{lutz} which contains 7 vertices and 14 facets:
	\begin{lstlisting}
		In[6]:= torus = SimplicialComplex["Torus"];
		In[7]:= torus["VertexCount"]
		Out[7]= 7
		In[8]:= torus["FacetCount"]
		Out[8]= 14
	\end{lstlisting}
	The torus can also be computed as the product $S^1\times S^1$ however this triangulation may not be minimal:
	\begin{lstlisting}
		In[9]:= SimplicialProduct[s1, s1]["FacetCount"]
		Out[9]= 18
	\end{lstlisting}
	
	The paclet supports creating new complexes from old, by taking joins or products of existing complexes (as seen above too). It also supports constructing the cone (join with a point) or suspension (join with two distinct points) on a simplicial complex. For example, cone of the $n$-dimensional simplex is the $n+1$-dimensional simplex:
	\begin{lstlisting}
		In[9]:=  n = SimplicialComplex[{"Simplex", 5}];
		In[10]:= np1 = SimplicialComplex[{"Simplex", 6}];
		In[11]:= cone = SimplicialCone[n];
		In[12]:= SimplicialIsomorphicQ[np1, cone]
		Out[12]= True
	\end{lstlisting}
	
	\subsection{Homology and Betti numbers}
	
	The primary feature of the paclet is the computation of simplicial homology groups. By default, computations are performed over $\mathbb{Z}$, which allows the package to detect topological torsion. We demonstrate this using the real projective plane, a non-orientable surface:
	\begin{lstlisting}
		In[14]:= rp2 = SimplicialComplex["RealProjectivePlane"];
		In[15]:= HomologyGroup[rp2]
		Out[15]= <|0 -> $\mathbb{Z}$, 1 -> $\mathbb{Z}_2$, 2 -> 0|>
	\end{lstlisting}
	We can also query the Betti numbers directly, which return the rank of the free abelian component. Furthermore, the paclet supports computations over rationals or finite fields. By the Universal Coefficient Theorem, computing the homology of $\mathbb{RP}^2$ over $\mathbb{F}_2$ alters the Betti numbers:
	\begin{lstlisting}
		In[16]:= BettiNumber[rp2] (* over integers by default *)
		Out[16]= <|0 -> 1, 1 -> 0, 2 -> 0|>
		In[17]:= BettiNumber[rp2, "Coefficients" -> Rationals]
		Out[17]= <|0 -> 1, 1 -> 0, 2 -> 0|>
		In[18]:= BettiNumber[rp2, "Coefficients" -> 2] (* over F_2 *) 
		Out[18]= <|0 -> 1, 1 -> 1, 2 -> 1|> 
	\end{lstlisting}
	
	\subsection{Large scale computation}
	
	To demonstrate the performance of the paclet below I evaluate the homology of all 180 unlabeled abstract simplicial complexes on up to 5 vertices:
	\begin{lstlisting}
		In[19]:= Length[facets]
		Out[19]= 180
		In[20]:= complexes = SimplicialComplex /@ facets;
		In[21]:= BarChart[Counts[#["EulerCharacteristic"] & /@ complexes], (* options *)]
		Out[21]= $\text{Shown in figure \ref{bar}.}$
	\end{lstlisting}
	The most common homology group:
	\begin{lstlisting}
		In[22]:= TakeLargest[Counts[HomologyGroup /@ complexes], 1]
		Out[22]= <|<|0 -> $\mathbb{Z}$, 1 -> $\mathbb{Z}$, 2 -> 0|> -> 27|>
	\end{lstlisting}
	This shows that out of all unlabelled abstract simplicial complexes on up to 5 vertices, the homology of $S^1$ appears the most frequently.
	\begin{figure}[h]
		\centering
		\includegraphics[width=.6\linewidth]{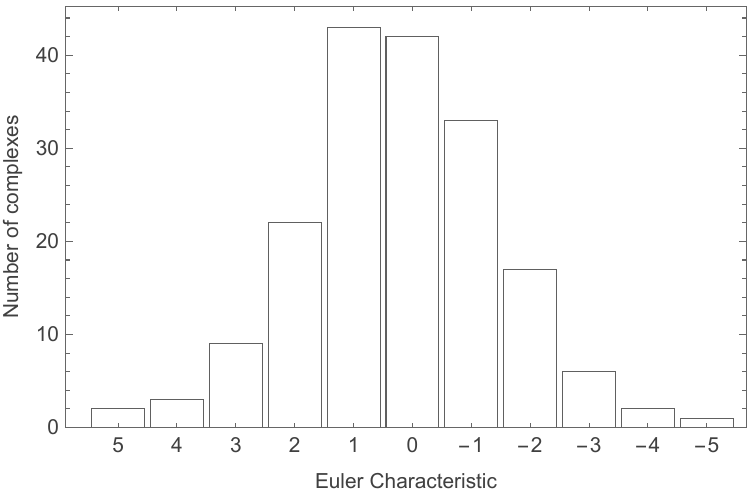}
		\caption{Out[21] showing the number of unlabelled abstract simplicial complexes having Euler characteristic varying from $-5$ to $5$.}
		\label{bar}
	\end{figure}
	
	To provide a representative indication of performance, table~\ref{bench} reports the execution times of several core operations on a collection of standard simplicial complexes of various ``sizes.'' As expected, the computational cost increases with the size of the underlying simplicial complex, while Betti number computations remain considerably faster than full homology computations since they avoid determining the torsion subgroup. The results demonstrate that the package is capable of efficiently handling a variety of complexes commonly encountered in algebraic topology.
	
	\begin{table}[ht]
		\centering
		\caption{Execution times for representative simplicial complexes. Environment: HP 15s, AMD Athlon Silver processor, 16 GB RAM, 512 GB SSD, Windows 11, Mathematica 15.0.}
		\label{bench}
		\begin{tabular}{>{\raggedright\arraybackslash}p{3cm}lrrr}
			\toprule
			Complex & $f$-vector & Homology & Betti & Aut. grp\\
			\midrule
			$S^1$ & $(1,3,3)$ & 0.0024 & 0.0020 & 0.0005 \\
			Torus & $(1,7,21,14)$ & 0.0109 & 0.0052 & 0.0008 \\
			Klein bottle & $(1,8,24,16)$ & 0.0100 & 0.0051 & 0.0007 \\
			$\mathbb{RP}^2$ & $(1,6,15,10)$ & 0.0074 & 0.0041 & 0.0007 \\
			$\mathbb{RP}^3$ & $(1,11,51,80,40)$ & 0.1904 & 0.0184 & 0.0017 \\
			$\mathbb{RP}^4$ & $(1,16,120,330,375,150)$ & 1.9004 & 0.3467 & 0.0069 \\
			Poincar\'e homology sphere & $(1,16,106,180,90)$ & 0.4557 & 0.0475 & 0.0019 \\
			\bottomrule
		\end{tabular}
	\end{table}
	
	\section{Conclusion}
	\label{con}
	
	This paper presented SimplicialHomology, an open-source Wolfram Language paclet for constructing abstract simplicial complexes and computing their homological invariants. The package provides functionality for computing homology groups over the integers, rational numbers, and finite fields, together with related invariants such as Betti numbers and Euler characteristics. In addition, it supports common constructions and operations on simplicial complexes, including joins, cones, suspensions, stars, links, connected components, and automorphism groups. The implementation employs sparse boundary matrices and appropriate caching to improve computational efficiency while remaining fully integrated with the symbolic features of the Wolfram Language. The examples and benchmarks demonstrate that the package is capable of efficiently handling a range of simplicial complexes encountered in teaching and research. Future work includes support for simplicial maps and induced homomorphisms, formalisation of chain complexes, and cohomology.
	
	\appendix
	\section{Availability}
	\label{a}
	The package is available on the Wolfram's official paclet repository at the link \url{https://resources.wolframcloud.com/PacletRepository/resources/Taggar/SimplicialHomology/} from where it can be installed in any local or cloud session of Mathematica. The source code of the package is available in its GitHub repository at \url{https://github.com/zplus11/SimplicialHomology/}. An archive of the package is available at \url{https://doi.org/10.5281/zenodo.21309175}.
	
	\bibliographystyle{plain}
	\bibliography{references}

\end{document}